# Topological Recursion Relations by Localization

Eric Edward Katz

2 Oct 2003

## 1 Introduction

Let $\overline{\mathcal{M}}_{g,n}$ be the moduli space of stable genus $g$ curves with $n$ marked points. $\overline{\mathcal{M}}_{g,n}$ has boundary strata consisting of nodal curves parametrized by products of smaller-dimensional moduli spaces. The fundamental classes of these boundary strata may be linearly dependent in $A_*(\overline{\mathcal{M}}_{g,n})$. Relations among these boundary strata can be found by exploiting a localization trick. Let $\pi$ be the stabilization map

$$\pi : \overline{\mathcal{M}}_{g,n}(\mathbb{P}^1, d) \to \overline{\mathcal{M}}_{g,n}.$$

Note that $\overline{\mathcal{M}}_{g,n}(\mathbb{P}^1, d)$ has a $\mathbb{C}^*$-action induced from the $\mathbb{C}^*$-action on $\mathbb{P}^1$. $\overline{\mathcal{M}}_{g,n}$ has a trivial $\mathbb{C}^*$-action which makes $\pi$ an equivariant map. If we write $A_*^{\mathbb{C}^*}(pt) = \mathbb{C}[\hbar]$ then $A_*^{\mathbb{C}^*}(\overline{\mathcal{M}}_{g,n}) = A_*(\overline{\mathcal{M}}_{g,n})[\hbar]$.

Now, given $a \in A^*_{\mathbb{C}^*}(\overline{\mathcal{M}}_{g,n}(\mathbb{P}^1, d))$, we may compute $\pi_*([\overline{\mathcal{M}}_{g,n}(\mathbb{P}^1, d)]^{vir} \cap a)$ by the following localization formula

$$\pi_* \left([\overline{\mathcal{M}}_{g,n}(\mathbb{P}^1, d)]^{vir} \cap a\right) = \sum_F \pi_* i_{F*} \left(\frac{[F]^{vir} \cap i_F^* a}{e(\nu)}\right) \tag{1}$$

in $A_*^{\mathbb{C}^*}(\overline{\mathcal{M}}_{g,n}) \otimes \mathbb{C}[\hbar, \frac{1}{\hbar}]$ where the sum is over $\mathbb{C}^*$-fixed loci $F$ in $\overline{\mathcal{M}}_{g,n}(\mathbb{P}^1, d)$, $i_F : F \to \overline{\mathcal{M}}_{g,n}(\mathbb{P}^1, d)$ and $e(\nu)$ is the equivariant Euler class of the normal bundle to $F$ in $\overline{\mathcal{M}}_{g,n}(\mathbb{P}^1, d)$.

Note that the left hand side of (1) is a finite expansion in powers of $\hbar$ and involves no negative powers of $\hbar$. The right hand side of (1) involves $\frac{1}{e(\nu)}$ and so has negative powers of $\hbar$ occuring. This means that the coefficient of $\hbar^{-k}$ on the right hand side must be zero. Because the coefficient is the sum of cycle classes pushed forward from fixed loci, we obtain relations among the boundary strata in $A_*(\overline{\mathcal{M}}_{g,n})$.

In what follows, we will apply this trick to

$$\pi : \overline{\mathcal{M}}_{1,4}(\mathbb{P}^1, 2) \to \overline{\mathcal{M}}_{1,4}$$

and obtain the Getzler relation in $A_*(\overline{\mathcal{M}}_{1,4})$

The Getzler relation was first discovered by Getzler in [2] by bounding the dimension of $H_4(\overline{\mathcal{M}}_{1,4})$ and computing its intersection matrix. A geometric proof of this relation was given by Pandharipande in [8]. Pandharipande's proof makes heavy use of the branch morphism and can probably also be cast in the more modern language of the degeneration methods of [6]. This note can be thought of as an attempt to modify the proof of Pandharipande to use the technique of virtual localization.

I would like to mention that the main technique of this note is used in a number of papers of Faber and Pandharipande, most recently, implicitly in section 2.2.2 of [1].

I would like to extend thanks to Ravi Vakil for valuable discussions and to Rahul Pandharipande for the use of his figures.

# 2 Dimension Count

The purpose of this section is to state the virtual dimensions of the relevant moduli spaces.

$$\mathrm{vdim}_{\mathbb{C}}\overline{\mathcal{M}}_{g,n} = 3g - 3 + n$$
$$\mathrm{vdim}_{\mathbb{C}}\overline{\mathcal{M}}_{g,n}(\mathbb{P}^1, d) = 2g - 2 + 2d + n$$

Therefore, the virtual dimension of the fiber of the stabilization map $\pi$ is

$$\mathrm{vdim}_{\mathbb{C}}\overline{\mathcal{M}}_{g,n}(\mathbb{P}^1, d) - \mathrm{vdim}_{\mathbb{C}}\overline{\mathcal{M}}_{g,n} = (2g-2+2d+n) - (3g-3+n) = 1 - g + 2d$$

Hence if $a \in A^m_{\mathbb{C}^*}(\overline{\mathcal{M}}_{g,n}(\mathbb{P}^1, d))$, $\pi_*([\overline{\mathcal{M}}_{g,n}(\mathbb{P}^1, d)]^{vir} \cap a) \in A_{2g-2+2d+n-m}(\overline{\mathcal{M}}_{g,n})$.

Since $\hbar$ is of degree 1 in $A^{\mathbb{C}^*}_*(\overline{\mathcal{M}}_{g,n})$, the coefficient of $\hbar^{-k}$ in $\pi_*([\overline{\mathcal{M}}_{g,n}(\mathbb{P}^1, d)]^{vir} \cap a)$ is an element of $A_{2g-2+2d+n-m-k}(\overline{\mathcal{M}}_{g,n})$.

# 3 Symmetrized Strata in $\overline{\mathcal{M}}_{1,4}$

The Getzler relation [2] is a relation among two-dimensional strata on $\overline{\mathcal{M}}_{1,4}$. It is best written in terms of strata that have been symmetrized over the placement of marked points. We adopt Getzler's notation for symmetrized strata, so the symbol

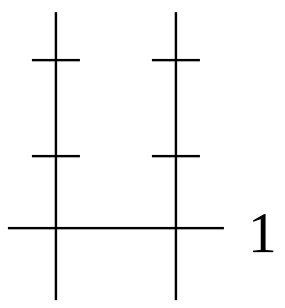

is shorthand for

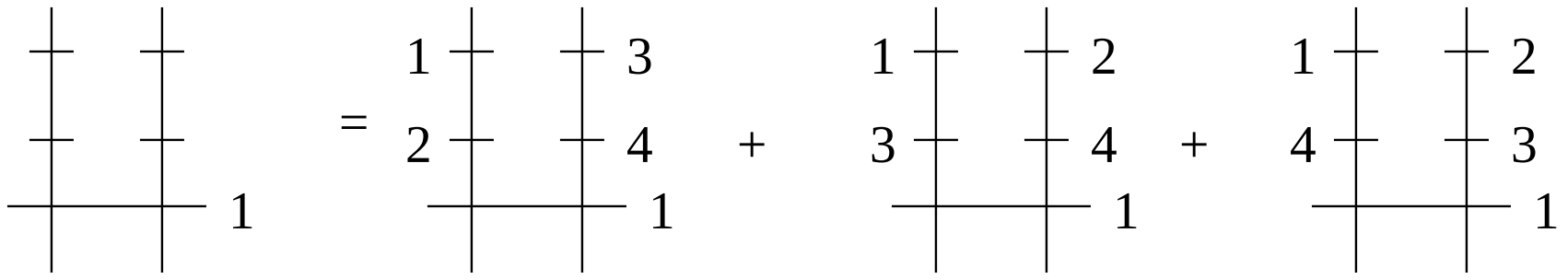

where the line labelled 1 is an elliptic curve, the unlabelled lines are rational curves and the bars represent marked points. In other words, given a stratum $S$ in $\overline{\mathcal{M}}_{1,4}$, we consider the symmetrized stratum

$$\mathrm{Sym}(S) = \frac{1}{\mathrm{Aut}(S)} \sum_{\sigma \in S_4} \sigma(S)$$

where $\sigma(S)$ is the action of $\sigma$ on $S$ permuting the marked points and $\mathrm{Aut}(S)$ is the number of permutations that leave $S$ fixed.

Getzler's notation for the set of $S_4$-invariant dimension 2 strata of $\overline{\mathcal{M}}_{1,4}$ is

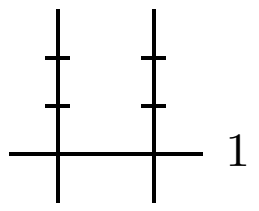

$\triangle_{2,2}$

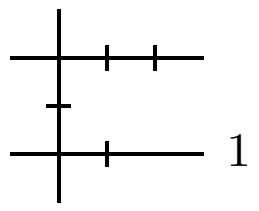

$\triangle_{2,3}$

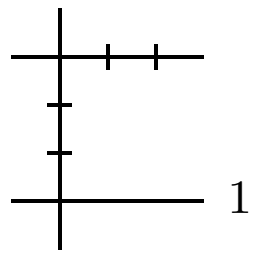

$\triangle_{2,4}$

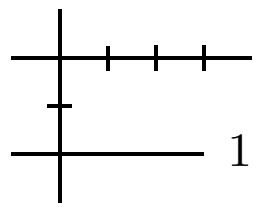

$\triangle_{3,4}$

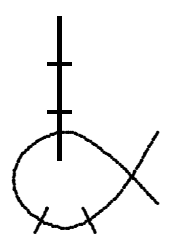

$\triangle_{0,2}$

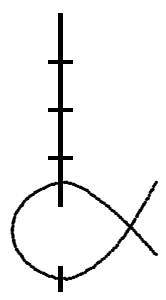

$\triangle_{0,3}$

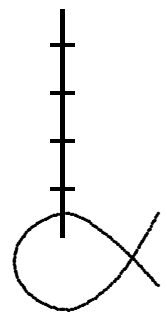

$\triangle_{0,4}$

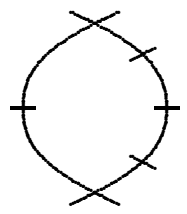

$\triangle_{a}$

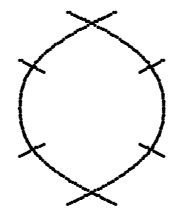

$\triangle_{b}$

Note tht these figures are taken from the appendix of [8].

# 4 The Getzler Relation

Using the above notation, we may write Getzler's relation as

$$12 \triangle_{2,2} - 4 \triangle_{2,3} - 2 \triangle_{2,4} + 6 \triangle_{3,4} + \triangle_{0,3} + \triangle_{0,4} - 2\triangle_{\beta} = 0. \qquad (2)$$

Note that this formula agrees with that in [2], not that in [8]. The two instances differ in that in [2], the symbols $\triangle$ refer to the fundamental class of the stratum considered as a stack, i.e. the usual fundamental class divided by the order of the automorphism group of the generic point while in [8], the same symbols refer to the usual fundamental class. We have found it convenient to use the fundamental classes of the stack throughout.

The Getzler relation appear as the coefficient of $\hbar^{-2}$ in $\pi_*([\overline{\mathcal{M}}_{1,4}(\mathbb{P}^1, 2)]^{vir} \cap a)$ for a particular choice of $a \in A^4_{\mathbb{C}^*}((\mathbb{P}^1, 2))$.

# 5 Group Actions and Linearizations

The purpose of this section is to review the standard conventions for the group action on $\mathbb{P}^1$. See [3] for more details. We define a $\mathbb{C}^*$ action on $\mathbb{P}^1$ by

$$\lambda \cdot [Z_0 : Z_1] = [\lambda Z_0 : Z_1]$$

We identify $\mathbb{C}$ with $\mathbb{P}^1 - [1 : 0]$ and refer to $[0 : 1]$ as 0 and $[1 : 0]$ as $\infty$.

We will pull back point classes from $\mathbb{P}^1$ by the evaluation map. To perform the localization computation, we need to pick equivariant extensions of the point class. These equivariant extensions will arise as $c_1(\mathcal{O}(1))$ for appropriate choices of linearization of $\mathcal{O}(1)$

If $\mathcal{O}(-1)$ is the line bundle whose fiber over $[Z_0 : Z_1]$ are all complex multiples of $(Z_0, Z_1)$, we define the linearization of $\mathcal{O}(-1)$ of weight $b$ to be

$$\lambda \cdot (Z_0, Z_1) = (\lambda^{b+1} Z_0, \lambda^b Z_1)$$

If $i_0 : 0 \to \mathbb{P}^1$ and $i_\infty : \infty \to \mathbb{P}^1$ are the usual inclusions then ${i_0}^*\mathcal{O}(-1)$ is the vector space $\mathbb{C}$ with the action

$$\lambda \cdot z = \lambda^{b+1} z$$

and ${i_\infty}^*\mathcal{O}(-1)$ is the vector space $\mathbb{C}$ with the action

$$\lambda \cdot z = \lambda^b z.$$

Let $\mathcal{O}(1)$ be the equivariant line-bundle on $\mathbb{P}^1$ dual to the one defined above. For $b = 0$, $c_1(\mathcal{O}(1)) = i_{0*}[P]$ where $P$ is Poincare-dual to the fundamental class of a point. Likewise, for $b = 1$, $c_1(\mathcal{O}(1)) = i_{\infty *}[P]$. Set $[0] = i_{\infty *}[P]$ and $[\infty] = i_{\infty *}[P]$. Note that ${i_0}^*[\infty] = 0$ and ${i_\infty}^*[0] = 0$.

The class $a \in A^4_{\mathbb{C}^*}(\overline{\mathcal{M}}_{1,4}(\mathbb{P}^1, 2))$ will be an appropriately symmetrized version of

$$a = ev_1^*[0] ev_2^*[0] ev_3^*[\infty] ev_4^*[\infty]$$

# 6 Fixed Loci

We will be summing over fixed loci of the $\mathbb{C}^*$-action on $\overline{\mathcal{M}}_{g,n}(\mathbb{P}^1, d)$. These loci are given by genus $g$ stable maps with $n$ marked points whose irreducible components are as follows:

(i) Curves of any genus contracted to 0

(ii) Curves of any genus contracted to $\infty$

(iii) Rational curves mapped into $\mathbb{P}^1$, ramified only over 0 and $\infty$, called trivial components

The marked points on the fixed loci are mapped to 0 or $\infty$.

An example of a fixed locus in $\overline{\mathcal{M}}_{1,4}(\mathbb{P}^1, 2)$ is

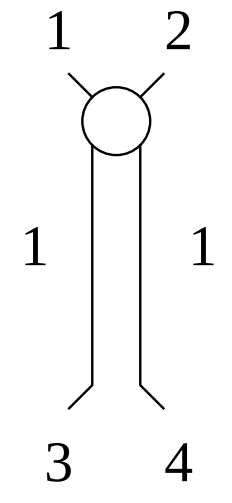

where the circle at top represents an elliptic curve contracted to 0, the 1's at the side are degree 1 trivial components, 1 and 2 are marked points on the elliptic curve, and the 3 and 4 marked points on the trivial components.

# 7 Evaluation of Classes on Fixed Loci

Following [5] or [3] we will compute

$$\sum_F \pi_* i_{F*} \left( \frac{[F]^{vir} \cap i_F^*(ev_1^*[0] ev_2^*[0] ev_3^*[\infty] ev_4^*[\infty])}{e(\nu)} \right)$$

by parametrizing fixed loci by a product of moduli spaces. For example, the above fixed locus is parametrized by $\overline{\mathcal{M}}_{1,4}$. More formally, we define a map

$$I : \overline{\mathcal{M}}_{1,4} \to \overline{\mathcal{M}}_{1,4}(\mathbb{P}^1, 2)$$

that attaches degree one trivial components to the third and fourth marked marked points and then places the marked points 3 and 4 on the trivial components.

So we note that $[F]^{vir} = \frac{1}{\deg(I)} I_*([\overline{\mathcal{M}}_{1,4}]^{vir})$ where $\deg(I)$ is the degree of the map $I$ considered as a map of stacks which in more conventional language is the number of automorphisms of $F$ divided by the number of automorphisms of $\overline{\mathcal{M}}_{1,4}$. In this case, $\deg(I)$ is equal to 1 but in more general situations, we may have more complicated fixed loci which are products of moduli spaces and whose automorphism groups are non-trivial. See [5] or [3] for more explanation. See also [4] for a very systematic description of fixed locii of $\overline{\mathcal{M}}_{g,n}(\mathbb{P}^1, d)$.

In general, we will have a product of moduli stacks $L$ parametrizing a fixed locus by a map $I : L \to F$. This gives us the following commutative diagram

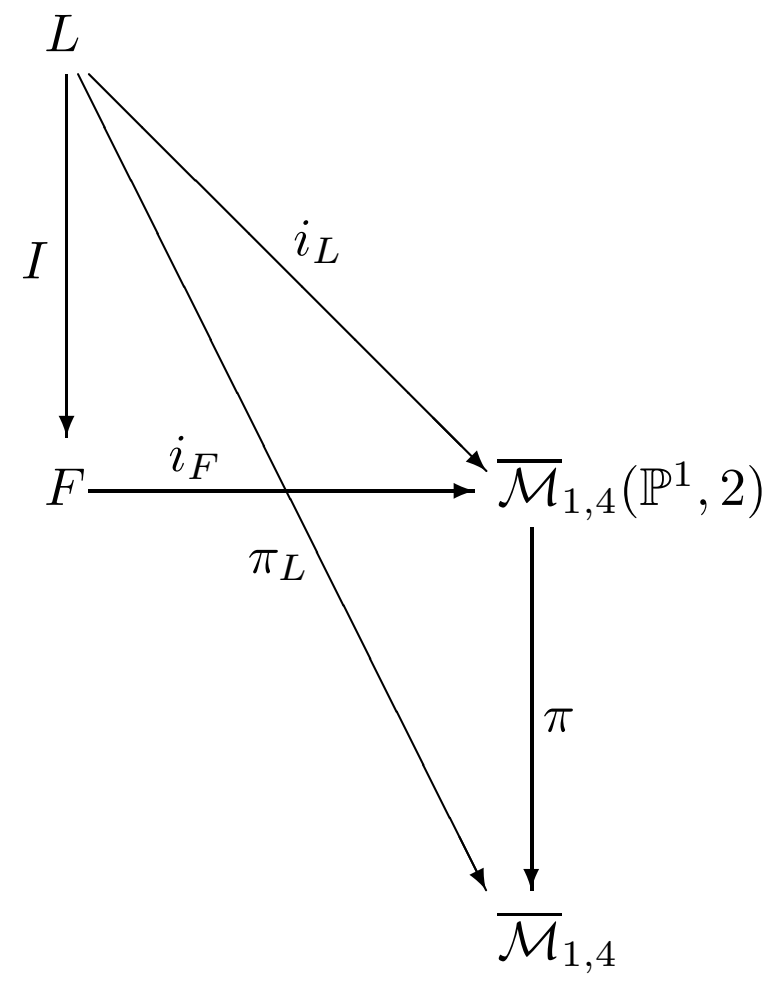

Therefore,

$$\pi_* i_{F*}\left(\frac{[F]^{vir}\cap i_F^* a}{e(\nu)}\right) = \frac{1}{\deg(i_L)}\pi_{L*}\left(\frac{[L]^{vir}\cap i_L^* a}{I^* e(\nu)}\right).$$

Now, note that $\pi_L$ parametrizes a stratum $S$ in $\overline{\mathcal{M}}_{1,4}$. Since we are using fundamental classes for strata considered as stacks,

$$[S]^{vir} = \frac{1}{\deg(\pi_L)}\pi_{L*}([L]^{vir}).$$

In what follows, we will abuse notation and identify $[S]^{vir}$ with the cycle class on $L$ which pushes forward by $\pi_{L*}$ to the legitimate $[S]^{vir}$. It follows that

$$\pi_* i_{F*}\left(\frac{[F]^{vir}\cap i_F^* a}{e(\nu)}\right) = \frac{\deg(\pi_L)}{\deg(i_L)}\pi_{L*}\left([S]^{vir}\cap \frac{i_L^* a}{I^* e(\nu)}\right).$$

Now, we notice that $i_L^*(ev_i^*[0])$ and $i_L^*(ev_i^*[\infty])$ are of pure weight. In fact, $i_L^*(ev_1^*[0]ev_2^*[0]ev_3^*[\infty]ev_4^*[\infty]) = \delta_F \hbar^4$ where $\delta_F = 1$ if marked points $1, 2$ are mapped to 0 and $3, 4$ are mapped to $\infty$ and $\delta_F = 0$ otherwise.

Therefore, we must evaluate

$$\sum_F \delta_F \frac{\deg(\pi_L)}{\deg(i_L)}\pi_{L*}\left(\frac{\hbar^4 [S]^{vir}}{I^* e(\nu)}\right). \qquad (3)$$

Since we are looking at the coefficient of $\hbar^{-2}$, in (3), we can consider the constant coefficient of

$$\sum_F \delta_F \frac{\deg(\pi_L)}{\deg(i_L)} \pi_{L*} \left( \frac{\hbar^6 [S]^{vir}}{I^* e(\nu)} \right). \tag{4}$$

# 8 Symmetrization of Marked Points

In (4), the cycles classes in each summand is independent of the labelling of marked points. Therefore,

$$\begin{aligned} & \sum_{\sigma\in S_4} \sum_F \pi_* i_{F*} \frac{[F]^{vir} \cap i_F^*(ev^*_{\sigma(1)}[0] ev^*_{\sigma(2)}[0] ev^*_{\sigma(3)}[\infty] ev^*_{\sigma(4)}[\infty])}{e(\nu)} \\ = & \sum_{\sigma\in S_4} \sum_{\sigma(F)} \delta_{\sigma(F)} \frac{\deg(\pi_L)}{\deg(i_L)} \pi_{L*} \left( \frac{\hbar^6 [S]^{vir}}{I^* e(\nu)} \right) \end{aligned}$$

where $\sigma$ acts on the right hand side by permuting marked points.

For ease of notation, we can replace fixed locii in $\overline{\mathcal{M}}_{1,4}(\mathbb{P}^1, 2)$ labelled with marked points by unlabelled fixed locii. In other words, we consider

$$\mathrm{Sym}(F) = \sum_{\sigma \in S_4} \sigma(F)$$

Note that our convention for symmetrized strata in $\overline{\mathcal{M}}_{1,4}(\mathbb{P}^1, 2)$ differs from that in $\overline{\mathcal{M}}_{1,4}$ and that

$$\pi_*(\mathrm{Sym}(F)) = \mathrm{Aut}(\pi_*(F))\ \mathrm{Sym}(\pi_* F)).$$

The unlabelled locus

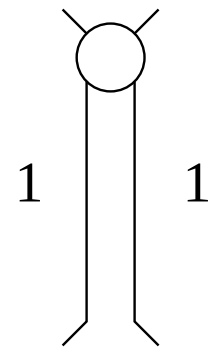

represents

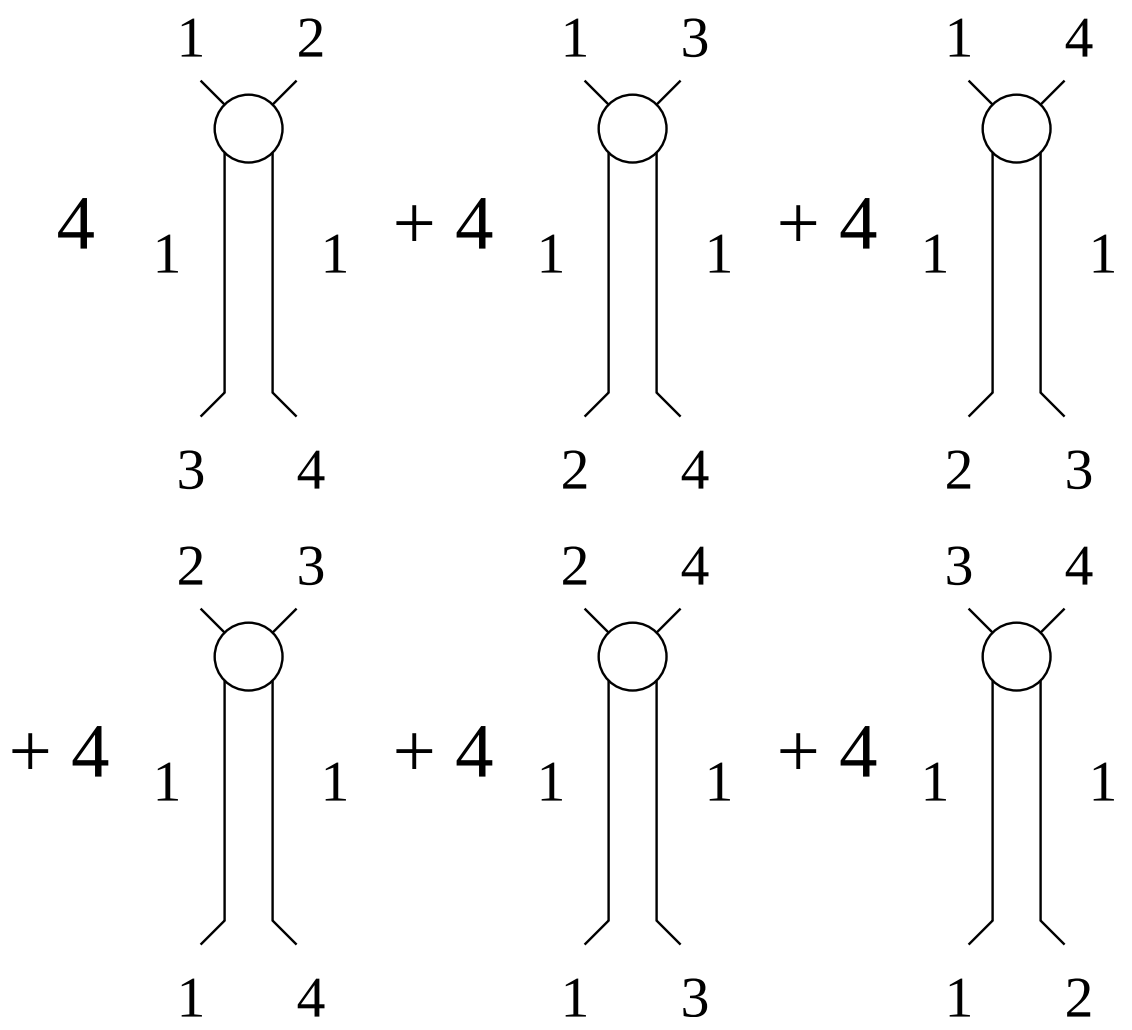

we consider, the shorthand of where we pick up a factor of 4 when we pushforward by $\pi_*$.

Now, instead of summing over labelled fixed locii, we sum over their unlabelled counterparts, i.e.

$$\sum_{\sigma\in S_4}\sum_{\sigma(F)}\delta_{\sigma(F)}\frac{\deg(\pi_L)}{\deg(i_L)}\pi_{L*}\left(\frac{\hbar^6[L]^{vir}}{I^*e(\nu)}\right)$$
$$= \sum_{F_{\mathrm{un}}}\frac{\deg(\pi_L)}{\deg(i_L)}N_{F_{\mathrm{un}}}Aut(\pi_*L)\pi_{L*}\left(\frac{\hbar^6[S]^{vir}}{I^*e(\nu)}\right)$$

where $N_{F_{\mathrm{un}}}$ is the number of labelled fixed with $\delta_F = 1$ locii have $F_{\mathrm{un}}$ as their unlabelled counterpart. More formally, this is

$$N_{F_{\mathrm{un}}} = \sum_{\mathrm{un}(F)=F_{\mathrm{un}}}\delta_F \tag{5}$$

where $\mathrm{un}(F)$ denotes the unlabelled fixed locus corresponding to $F$.

For example,

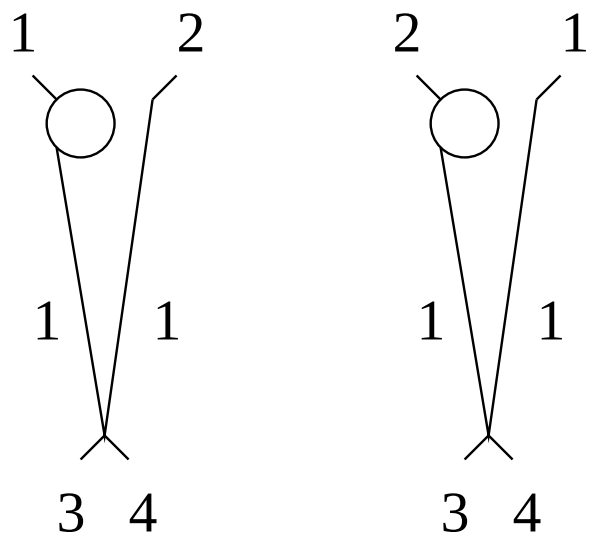

which are distinct fixed locii in $\overline{\mathcal{M}}_{1,4}(\mathbb{P}^1, 2)$ both have as their unlabelled counterpart.

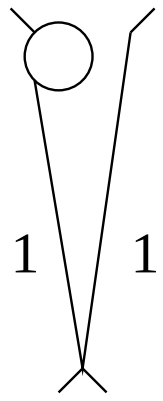

# 9 Contributing Fixed Locii

The unlabelled fixed locii in $\overline{\mathcal{M}}_{1,4}(\mathbb{P}^1, 2)$ are

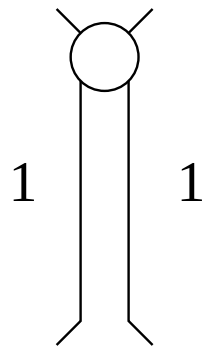

$\prod_{1,1}$

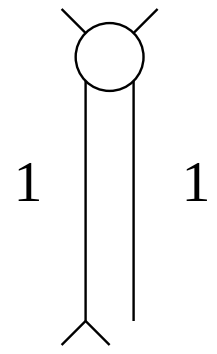

$\prod_{2,0}$

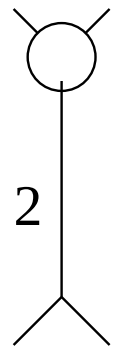

$\top$

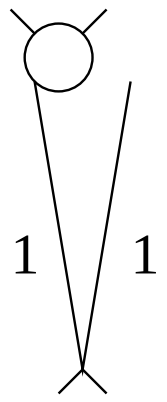

$\mathsf{V}_{2,0}$

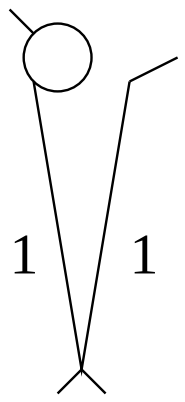

$\mathsf{V}_{1,1}$

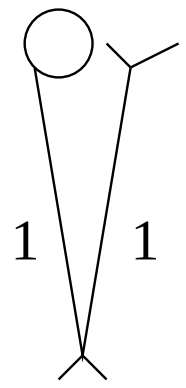

$\mathsf{V}_{0,2}$

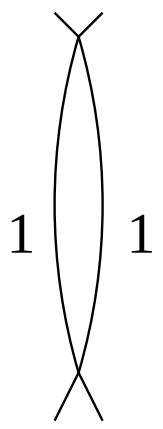

$O$

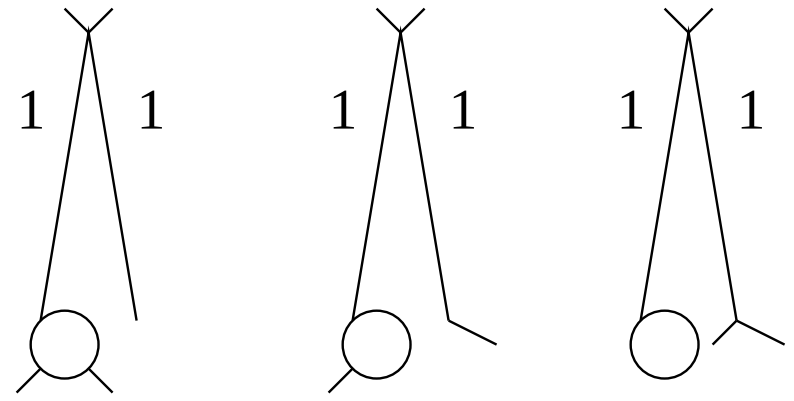

$\bigwedge_{2,0}$ $\bigwedge_{1,1}$ $\bigwedge_{0,2}$

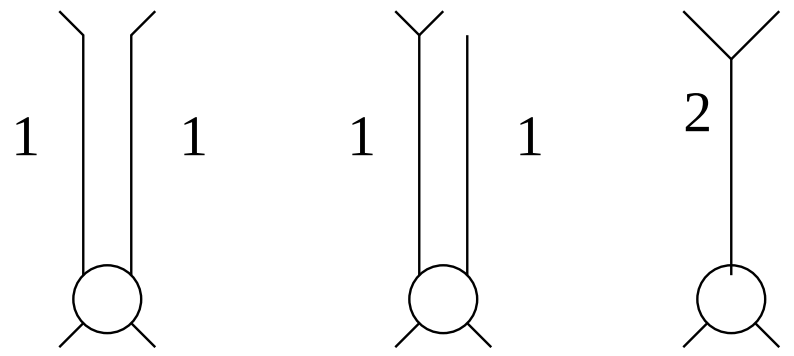

$\coprod_{1,1}$ $\coprod_{2,0}$ $\perp$

# 10 Example Computation

Determining the contribution of each fixed locus is straight-forward but laborious. To give the reader an idea of the computation, we will find the contribution of $\bigvee_{1,1}$

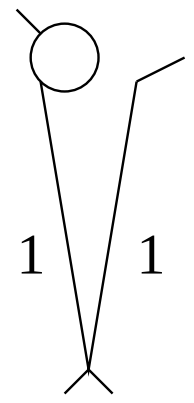

Note that this fixed locus is paramterized by $\overline{\mathcal{M}}_{1,2} \times \overline{\mathcal{M}}_{0,4}$ where the elliptic curve is contracted to 0 and the rational curve is contracted to $\infty$. This locus has $N_F = 2$. We must determine the normal bundle to this fixed locus. This is done in detail in [5] and [3] so we will only give a very colloquial explanation here. We can think of the normal bundle as being given by

$$e(\nu) = \frac{e(\text{Node Resolutions})e(\text{Moving})}{e(\text{Hodge})e(\text{Structure})e(\text{Automorphisms})} \tag{6}$$

where Node Resolutions denote the bundle of deformations given by resolving nodes, Moving, the bundle of deformations of the map from the trivial components, Hodge, the appropriate Hodge bundle, Structure, the bundle of equations on deforming the map that force the map to maintain its dual graph, and Automorphisms, the bundle of automorphisms of the curve. Here,

$$e(\text{Node Resolutions}) = (-\psi_{\text{ell}} + \hbar)(-\psi_{\text{rat},1} - \hbar)(-\psi_{\text{rat},2} - \hbar)$$

where $\psi_{\text{ell}}$ is the $\psi$-class at the node on the elliptic curve and $\psi_{\text{rat},i}$ are the $\psi$-classes at the nodes on the rational curve.

$$e(\text{Moving}) = ((-\hbar)(\hbar))^2$$

$$e(\text{Hodge}) = \hbar - \lambda_1$$

where $\lambda_1$ is the first Hodge class on the moduli of elliptic curves.

$$e(\text{Structure}) = -\hbar$$

which comes from the equation that forces the two trivial components to meet at the collapsed rational curve.

$$e(\text{Automorphisms}) = 1$$

Therefore, we get

$$e(\nu) = \frac{(\hbar + \psi_{\text{rat},1})(\hbar + \psi_{\text{rat},2})(\hbar - \psi_{\text{ell}})(\hbar^4)}{(\hbar - \lambda_1)(-\hbar)}.$$

On $L = \overline{\mathcal{M}}_{1,2} \times \overline{\mathcal{M}}_{0,4}$, the contribution is given by

$$\frac{\hbar^6}{e(\nu)} = -\frac{\hbar^6(\hbar - \lambda_1)\hbar}{(\hbar + \psi_{\text{rat},1})(\hbar + \psi_{rat,2})(\hbar - \psi_{\text{ell}})\hbar^4} = -\frac{\hbar(1 - \frac{\lambda_1}{\hbar})}{(1 + \frac{\psi_{\text{rat},1}}{\hbar})(1 + \frac{\psi_{\text{rat},2}}{\hbar})(1 - \frac{\psi_{\text{ell}}}{\hbar})}.$$

The coefficient of $\hbar^0$ in the above is

$$\lambda_1 + \psi_{\text{rat},1} + \psi_{\text{rat},2} - \psi_{\text{ell}} \tag{7}$$

By standard string equation arguments, $\psi_{\text{rat},1}$ and $\psi_{rat,2}$ each force the rational curve to become a two component nodal curve with two special points on each components. Therefore, $\psi_{rat,1} + \psi_{rat,2}$ evaluated on the above fixed locus is

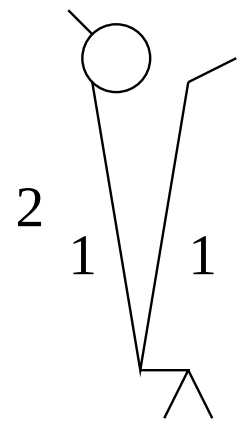

where the short horizontal line is a node. In other words, this is the stratum parametrized by $\overline{\mathcal{M}}_{1,2}$ which consists of a one-pointed elliptic curve with contracted to 0, a trivial component joining the elliptic curve to a two component nodal rational curve. Attached to one irreducible component of the nodal rational curve, there are two trivial component. Attached to the other component are two marked points mapping to $\infty$. And on the preimage of the trivial component disjoint from the elliptic curve is a marked point that maps to 0. This pushes forward to $2\triangle_{2,3}$. Since $\mathrm{Aut}(\pi_* L) = 2$ and $N_{F_{\mathrm{un}}} = 2$, we get a contribution of $8\triangle_{2,3}$.

On $\overline{\mathcal{M}}_{1,1}$, $\psi_{\mathrm{ell}} = \lambda_1$, so on $\overline{\mathcal{M}}_{1,2}$, $\psi_{\mathrm{ell}}$ and $\lambda_1$ differ only by a correction given by the string equation. It follows that $\psi_{\mathrm{ell}} - \lambda_1$ on $\overline{\mathcal{M}}_{1,2}$ is given by the 1-dimensional locus where the two marked points are on a rational component together. Therefore the contribution of $\lambda_1 - \psi_{\mathrm{ell}}$ is

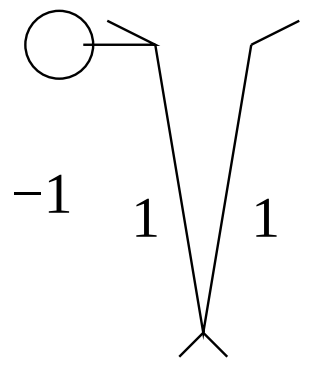

This is parametrized by $\overline{\mathcal{M}}_{1,1} \times \overline{\mathcal{M}}_{0,4}$ and pushes foward to $-\triangle_{3,4}$ with $\mathrm{Aut}(\pi_* L) = 6$ and $N_{F_{\mathrm{un}}} = 2$, so we get a contribution of $-12\triangle_{3,4}$

# 11 Contributions

The result of the calculations for all the fixed locii is

$$\bigvee_{0,2} : -4 \triangle_{2,4}$$

$$\bigvee_{1,1} : -12 \triangle_{3,4} +8 \triangle_{2,3}$$

$$\bigvee_{2,0} : 8 \triangle_{2,3} +8 \triangle_{2,4}$$

$$\top : -\frac{8}{3} \triangle_{0,2} -64 \triangle_{2,3} -64 \triangle_{2,4}$$

$$\prod_{2,0} : 8 \triangle_{2,2} +20 \triangle_{2,3} +24 \triangle_{2,4} + \triangle_{0,2}$$

$$\prod_{1,1} : 16 \triangle_{2,2} +20 \triangle_{2,3} +32 \triangle_{2,4} +24 \triangle_{3,4} +\frac{5}{3} \triangle_{0,2} +2 \triangle_{0,3} +2 \triangle_{0,4}$$

$$\bigwedge_{0,2} : -4 \triangle_{2,4}$$

$$\bigwedge_{1,1} : -12 \triangle_{3,4} +8 \triangle_{2,3}$$

$$\bigwedge_{2,0} : 8 \triangle_{2,3} +8 \triangle_{2,4}$$

$$\perp : -\frac{8}{3} \triangle_{0,2} -64 \triangle_{2,3} -64 \triangle_{2,4}$$

$$\coprod_{2,0} : 8 \triangle_{2,2} +20 \triangle_{2,3} +24 \triangle_{2,4} + \triangle_{0,2}$$

$$\coprod_{1,1} : 16 \triangle_{2,2} +20 \triangle_{2,3} +32 \triangle_{2,4} +24 \triangle_{3,4} +\frac{5}{3} \triangle_{0,2} +2 \triangle_{0,3} +2 \triangle_{0,4}$$

$$O : -8\triangle_b$$

Their sum is

$$4(12 \triangle_{2,2} -4 \triangle_{2,3} -2 \triangle_{2,4} +6 \triangle_{3,4} + \triangle_{0,3} + \triangle_{0,4} -2\triangle_b) = 0.$$

This is the Getzler relation.

## 12 Conclusion

The above method can be generalized and systematized in many directions. An obvious question is to apply this method to derive the Belorousski-Pandharipande relation on $\overline{\mathcal{M}}_{2,3}$. Another question are to systematize this method so that the computation can be done by a computer or can be expressed in terms of a Feynman expansion of an integral. A humbler version of the above is to prove that there are non-trivial relations by automating the computation of the leading terms (the symbol) of such relations. One also would like to study the structure of such relations. It is likely that much work will be done along these lines by using localization and degeneration methods in conjunction as in [1] and [7].

Department of Mathematics
Stanford University
eekatz@math.stanford.edu